\documentclass [11pt,a4paper]{article}
\usepackage{amsthm}
\usepackage{amsmath}
\usepackage{amssymb}
\usepackage{latexsym}
\usepackage{enumerate}

\newcommand{\Nbb}{\ensuremath{\mathbb{N}}}

\newcommand{\Rbb}{\ensuremath{\mathbb{R}}}
\newcommand{\Hbb}{\ensuremath{\mathbb{H}}}
\renewcommand{\Bbb}{\ensuremath{\mathbb{B}}}

\newcommand{\dcal}{\ensuremath{\mathcal{D}}}
\newcommand{\Hn}{\Hbb^n}
\newcommand{\Bn}{\Bbb^n}

\newcommand{\aster}{*\,}

\newcommand{\const}{\operatorname{const}}
\newcommand{\Iso}{\operatorname{Iso}}
\newcommand{\Id}{\operatorname{Id}}
\newcommand{\Tr}{\operatorname{Tr}}

\newcommand{\dsum}{\displaystyle\sum}
\newcommand{\rvec}{\stackrel{\to}{r}}

\newcommand{\Ndash}{\nobreakdash--}
\newcommand{\guio}[1]{\nobreakdash-\hspace{0pt}#1}

\newcommand{\forms}{\guio{forms}}
\newcommand{\form}{\guio{form}}

\newcommand{\invariant}{\guio{invariant}}

\newcommand{\admissible}{\guio{admissible}}

\newtheorem {theorem} {Theorem}[section]
\newtheorem {proposition} [theorem]{Proposition}

\newtheorem {lemma}  [theorem]{Lemma}

\newtheorem{MainTheorem}{Theorem}

\theoremstyle{definition}

\title{$L^p$-estimates for Riesz
transforms on forms in the Poincar\'{e} space}
\author{Joaquim Bruna\thanks{supported by projects BFM2002-04072-C02-02 and
2001SGR00172}\\Departament de Matem\`atiques\\ Universitat
Aut\`onoma de Barcelona}
\date{}

\begin{document}
\maketitle
\begin{abstract}
Using hyperbolic form convolution with doubly isometry-invariant
kernels, the explicit expression of the inverse of the de Rham
laplacian $\Delta$ acting on $m$-forms in the Poincar\'{e} space $\Hn$
is found. Also, by means of some estimates for hyperbolic singular
integrals, $L^p$-estimates for the Riesz transforms
$\nabla^i\Delta^{-1}, i\leq 2$, in a range of $p$ depending on
$m,n$ are obtained. Finally, using these, it is shown that
$\Delta$ defines topological isomorphisms in a scale of Sobolev
spaces $H^s_{m,p}(\Hn)$ in case $m\neq\frac{n\pm 1}{2}, \frac n2$.

{\bf Mathematics Subject Classification(2000)}: 53C21, 58J05,
58J50, 58J70

{\bf Key words}: Hodge-de Rham laplacian, Sobolev spaces, Riesz
transforms, hyperbolic form convolution
\end{abstract}

\section{Statement of results and preliminaries}\label{sec1}

{\bf 1.1.} The main object of study in this paper is the Hodge-de
Rham Laplacian $\Delta$ acting on $m$-forms in the Poincar\'{e}
hyperbolic space $(\Hn, g)$. The aim is to prove that $\Delta$
defines topological isomorphisms in a range $H^s_{m,p}(\Hn)$ of
Sobolev spaces of forms defined as follows. For \linebreak $0\leq
m\leq n,$ $1\leq p <\infty$ and $s\in\Nbb$, the Sobolev space
$H^s_{m,p}(\Hn)$ is the completion of the space $\dcal_m(\Hn)$ of
smooth $m$\forms\ with compact support with respect  the norm
\[
        \|\eta\|_{p,s} =\sum_{i=0}^s \|\nabla^{(i)}\eta\|_p\,.
\]
Here $\nabla^{(i)}$ means the $i$-th covariant differential of
$\eta$, and for a covariant tensor $\alpha$
\[
    \|\alpha\|_p=\left(\int_{\Hn}|\alpha(x)|^p\,d\mu(x)\right)^{\frac 1p}\,,
\]
$|\alpha|$ being the pointwise norm of $\alpha$ with respect the
metric $g$ and $d\mu$ the volume-invariant measure on $\Hn$ given
by $g$. The space $H^s_{m,p}(\Hn)$ can be alternatively defined in
terms of weak derivatives. The main result of this paper is:
\begin{MainTheorem}
$\Delta$ is a topological isomorphism from $H^{s+2}_{m,p}(\Hn)$ to
$H^s_{m,p}$ for $p\in (p_1,p_2)$ with
\[
   p_1=\frac{2(n-1)}{n-2+|n-2m|}\,\qquad\quad\frac{1}{p_1}+\frac{1}{p_2}=1
\]
in case $m\neq\frac{n\pm 1}{2},$ $\frac n2$.
\end{MainTheorem}

 In the exceptional case $m=\frac{n\pm
1}{2}$, $\Delta$ is one to one but is not a topological
isomorphism for any $p$. For this case we obtain as well some
weighted estimates. If $m=\frac n2$, $\Delta$ is known to have a
non-trivial kernel. Of course, Sobolev spaces $H^s_{m,p}$ can be
considered for non integer $s$ as well, and the same results holds
by interpolation.

Notice that the Hodge star operator $*$ establishes an isometry
from $H^s_{m,p}(\Hn)$ to $H^s_{n-m,p}(\Hn)$ which commutes with
$\Delta$, and this is why the range $(p_1, p_2)$ depends only on
$|n-2m|$. Notice too that the range $(p_1,p_2)$ always contains
$p=2$ in the non-critical case $|n-2m|>1$ and that for functions
($m=0$), the range of $p$ is $(1,\infty)$ (see comments below). We
point out that the range $(p_1,p_2)$ equals $|\frac 1p-\frac
12|<\frac{\sqrt \mu}{n-1}$, where $\mu$ denotes the greatest lower
bound for the spectrum of $\Delta$ in $H^0_{m,2}(\Hn)$, whose
value (\cite{D}) is $\mu=\frac{(n-1-2m)^2}{4}$ (for $m<\frac n2
$).

 For the Sobolev spaces for $p=2$, $H^s_{m,2}(\Hn)$, another proof of the theorem, based in energy
methods and valid for an arbitrary complete hyperbolic manifold,
is given in \cite{BG}. The motivation for the theorem, as with
\cite{BG}, comes from mathematical physics, where most operators
exhibit $\Delta$ as their principal part and results like the
above become essential to establish existence and uniqueness
theorems.

Our method of proof is simply to construct an explicit inverse $L$
for $\Delta$ on $\dcal_m(\Hn)$ and show that there is a gain of
two covariant derivatives
\[
        \|L\eta\|_{p,s+2}\leq \const \|\eta\|_{p,s}\,.
\]
Thus $L\eta$ plays the role of the classical Riesz transform in
the Euclidean setting. The most delicate part is of course
\[
    \|\nabla^{(2)}L\eta\|_p\leq \const\|\eta\|_p\,,\qquad
    p_1<p<p_2\,,\quad\eta\in\dcal_m(\Hn)\,.
\]

Riesz-type operators such as $\nabla\Delta^{-\frac 12},
\nabla^{(2)}\Delta^{-1}$ have extensively been studied in
different contexts, for the case of \emph{functions}. On symmetric
spaces, they are bounded in $L^p, 1<p<\infty$ and of weak type
$(1,1)$. This was shown in \cite{AL} for the first order ones in
some spaces, and later extended to all symmetric spaces in
\cite{A}. The $L^p$-boundedness holds as well for higher order
Riesz transforms in symmetric spaces, but not generally the weak
type $(1,1)$ estimate. In more general contexts, this has been
shown in \cite{L1}, \cite{L2}, \cite{L3}, among others. In case of
$m$-forms, $0<m<n$, as far as we know, there are much less known
results, and is for those that our result is new. In \cite{P1},
\cite{P2} some aspects of harmonic analysis of forms are
developed; in particular, the exact expression for the heat kernel
is given, and it is very likely that from it one can get as well
an explicit expression for $\Delta^{-1}$. Strictly speaking, to
prove the result an exact expression of $\Delta^{-1}$ is not
needed, it is enough having estimates for the resolvent both local
and  at infinity. In \cite{L3}, estimates of this kind are
obtained and applied to Sobolev-type inequalities for forms, and
they might work for this purpose too. However, we feel that our
approach, that we next describe, is more elementary and might be
interesting in itself.

The de Rham Laplacian $\Delta$ is invariant by all isometries
$\varphi$ of $\Hn$. These form a a group that we denote here by
$\Iso(\Hn)$. Denoting by $\varphi^*(\eta)(x)= \eta(\varphi(x))$
the pull-back of a form $\eta$ by $\varphi$, this means that
$\Delta$ and $\varphi^*$ commute, for all $\varphi\in\Iso(\Hn)$.
Therefore the inverse $L$ of $\Delta$ should commute too with
$\Iso(\Hn)$. Among all isometries of $\Hn$, the \emph{hyperbolic
translations} $\Tr(\Hn)$ constitute a (noncommutative) subgroup,
in one to one correspondence with $\Hn$ itself. In section~2 we do
some harmonic analysis for forms in $\Hn$ and introduce
\emph{hyperbolic convolution of forms} to describe all operators
acting on $m$\forms\  and commuting with $\Tr(\Hn)$. In a second
step (subsection~2.2) we characterize the hyperbolic convolution
kernels $k(x,y)$ corresponding to operators commuting with the
full group $\Iso(\Hn)$.

Once the general expression of an operator commuting with
$\Iso(\Hn)$ has been found, we look for our $L$ among these. This
corresponds to $L$ having a kernel $k(x,y)$ which is a
\emph{fundamental solution} of $\Delta$ in a certain sense, and
having the best decay of infinity. This kernel turns out to be
unique for $m\neq\frac{n\pm 1}{2}$, $\frac n2$, we call it the
\emph{Riesz kernel for $m$\forms\ in $\Hn$}, it is found in
subsection~3.1 and estimated in subsection 3.2. Section 4 is
devoted to the proof of the $L^p$-estimates. Here we use standard
techniques in real analysis (Haussdorf-Young inequalities, Schur's
lemma, etc.). For the second-order Riesz transform, to show its
boundedness in the specified range $(p_1,p_2)$ needs considering
some notion of ``hyperbolic singular integral''. There exist some
references dealing with this, e.g. \cite{L4}, \cite{I}, and giving
some criteria for $L^p$-boundedness that might apply; however, as
the singular integral arises locally we have found easier and more
elementary to treat it with the classical Euclidean
Calder\'{o}n-Zygmund theory as a local model, and patch it in a
suitable way to infinity.

\vspace{1cm}

\noindent {\bf 1.2.} We collect here several notations and known
facts about $\Hn$. We will use both the unit ball model $\Bn$ with
metric $g=4(1-|x|^2)^{-2}\sum_i dx^i\,dx^i$ and the half-space
model $\Rbb_+^n=\{x_n>0\}$ with metric $g=x_n^{-2}\sum_i
dx^i\,dx^i$. Both models are connected via the Cayley transform
$\psi\colon \Rbb_+^n\to\Bn$ given in coordinates by
\[
        y_i=\frac{2x_i}{\dsum_i^{n-1}x_i^2+(x_n+1)^2}\,,
        \quad i=1,\dots, n-1\,;\qquad\quad
        y_n=\frac{\dsum_1^n
        x_i^2-1}{\dsum_1^{n-1}x_i^2+(x_n+1)^2}\,.
\]
We denote by $e\in\Hn$ the point $(0,0,\dots,1)\in\Rbb_+^n$ or
$0\in\Bn$.

The metric $g$ defines a pointwise inner product
$(\alpha,\beta)(x)$ between forms at $x$, for every $x\in\Hn$, and
a volume measure $d\mu$. In the ball model $d\mu$ is written
\linebreak $d\mu(x)=2^n(1-|x|^2)^{-n} dx^1\dots dx^n$, and
$d\mu(x)=x_n^{-n}dx^1\dots dx^n$ in the half-space model. We
denote by $\langle\;,\;\rangle$ the pairing between forms that
makes $H_{m,2}^s(\Hn)$ a Hilbert space
\[
    \langle\alpha,\beta\rangle=\int_{\Hn}
    (\alpha,\beta)(x)\,d\mu(x)\,.
\]
We write $|\alpha|$ and $\|\alpha\|$ for the pointwise and global
norms, respectively, of the form $\alpha$. In terms of the Hodge
star operator $*$ the inner product can be written too
\[
        \langle\alpha,\beta\rangle=\int_{\Hn} \alpha\wedge
        \aster \beta\,.
\]

The group $\Tr(\Hn)$ of hyperbolic translations is in one to one
correspondence $x\mapsto T_x$ with $\Hn$ through the equation
$T_x(e)=x$. The equations of $z=T_xy$ are better described in the
half-space model by
\[
        z_i=x_ny_i+x_i\,,\quad i=1,\dots, n-1\,;\quad
        z_n=x_ny_n\,.
\]
It is easily checked that indeed $\Tr(\Hn)$ is a (non-commutative)
group. The inverse transformation of $T_x$ will be denoted $S_x$.
Another explicit isometry $\varphi_x$ mapping $e$ to $x$,
satisfying $\varphi_x^{-1}=\varphi_x$ is given in the ball model
by
\begin{equation}\label{eq1}
    \varphi_x(y)=\frac{(|x|^2-1)y+(|y|^2-2xy
    +1)x}{|x|^2|y|^2-2xy+1}\,.
\end{equation}
Since the isotropy group of $0$ is the orthogonal group $O(n)$,
the general expresion of $\varphi\in\Iso(\Hn)$ is
$\varphi=\varphi_x\circ U$, with $x=\varphi(0)$.

The hyperbolic (or geodesic) distance between $x,y\in\Hn$ is
written $d(x,y)$. We will rather use the \emph{pseudohyperbolic
distance} $r=r(x,y)$, related to $d$ by the formula $d(x,y)=2
\text{arctanh}\, r(x,y)$. The explicit expression of $r(x,y)^2$ in
the $\Rbb_+^n$ model and the $\Bn$ model is respectively
\begin{align}\label{eq2}
    r^2&=\frac{|x-y|^2}{|x-y|^2+4x_ny_n}\,,\quad
    x,y\in\Rbb_+^n\,,\nonumber\\[2mm]
    r^2&=|\varphi_x(y)|^2=\frac{|x-y|^2}{(1-|x|^2)(1-|y|^2)+|x-y|^2}\,,\qquad
    x,y\in\Bn\,.
\end{align}

Associated to the group of translations we have the basis of
orthonormal translation-invariant vector fields $X_i(x)=(T_x)_*
(X_i(e))$, such that $X_i(e)=\frac{\partial}{\partial x_i}$. They
satisfy $X_i(u\circ T_x)=(X_iu)\circ T_x$ for every smooth
function $u$. We will denote by $w^i(x)$ the dual basis of $X_i$,
which accordingly is orthonormal and tranlation invariant too:
$T_x^* w^i=w^i$. Their expression in the $\Rbb_+^n$ model is
simply
\[
    X_i(x)= x_n\frac{\partial}{\partial x_i}\,,\qquad
    w^i(x)=x_n^{-1}dx^i\,,\quad i=1,\dots,n\,.
\]
Because of their translation-invariance property, the $(X_i, w^i)$
are more suitable than the $(X_i,\eta^i)$ defined in the ball
model $\Bn$ by
\[
    Y_i(x)=\frac{(1-|x|^2)}{2}\frac{\partial}{\partial x_i}\,,
    \qquad \eta^i(x)=2(1-|x|^2)^{-1} dx^i\,.
\]
For an increasing multindex $I$ of length $|I|=m$ we write
$w^I=w^{i_1}\wedge w^{i_2}\wedge\dots\wedge w^{i_m}$, and
similarly $dx^I$ or $\eta^I$. The $\{w^I\}_I$ is an orthonormal
translation-invariant basis of $m$\forms.

Recall that the de Rham Laplacian is defined
$\Delta=d\delta+\delta d$, where $\delta$ is the adjoint of $d$
with respecto to $\langle\;,\;\rangle$. Although strictly speaking
not needed, the following expression of $\Delta$ in
$w^I$-coordinates will simplify the analysis at some points. If
$\alpha=\sum_I \alpha_Iw^I$, a computation shows that in case
$n\not\in J$
\begin{equation}\label{eq3.1}
        (\Delta \alpha)_J=\Delta\alpha_J + 2\sum_{k\in J}
        X_k\alpha_{Jk}-p(n-p-1)\alpha_J\,.\tag{3.1}
\end{equation}
Here $Jk$ means the multindex obtained replacing $k$ by $n$. In
case $n\in J$,
\begin{equation}\label{eq3.2}
    (\Delta\alpha)_J= \Delta\alpha_J - 2\sum_{l\not\in J} X_l \alpha_{lJ}
    - (1-p) (p-n)\alpha_J \tag{3.2}
\end{equation}
where $lJ$ means the multiindex obtained replacing $n$ by $l$. For
a function $f$
\[
    \Delta f = -\sum_{i=1}^n X_i^2 f +(n-1)X_n f\,.
\]
In the ball model, with usual coordinates,
\begin{equation}\label{eq3.3}
    \Delta f =-\frac 14 (1-|x|^2)^2 \sum_{i,j=1}^n \frac{\partial^2 f}{\partial x_i \partial x_j}
+(1-\frac n2) (1-|x|^2) \sum x_i \frac{\partial f}{\partial
x_i}\tag{3.3}
\end{equation}

\setcounter{equation}{3}

\section{Translation invariant and isometry invariant operators on forms}

\noindent {\bf 2.1.} We are interested in finding the general
expression of an operator acting on $m$\forms, and
isometry-invariant. In a first step we consider
\emph{translation-invariant operators} acting an $m$\forms; these
are described by what we might call \emph{hyperbolic convolution}
as follows. Let $k(x,y)$ be a double $m$\form\ in $x,y$ and define
\[
    (C_k\alpha)(x)=\int_{\Hn}\alpha(y)\wedge *_y k(x,y) =\langle
    \alpha, k(x,\cdot)\rangle\,,\qquad \alpha\in\dcal_m(\Hn)\,.
\]
If $T_z$ is a translation with inverse $S_z$
\begin{align*}
    C_k(T_z^*\alpha)(x)
        &= \int_{\Hn}(T_z^*\alpha)(y)\wedge *_yk(x,y)
            =\int_{\Hn}\alpha(T_z y)\wedge *_y k(x,y)\\[2mm]
        &\qquad\qquad\qquad =\int_{\Hn}\alpha(y)\wedge *_y k(x,S_z y)\\[3mm]
   T_z^*(C_k\alpha)(x)&= C_k\alpha(T_z
   x)=\int_{\Hn}\alpha(y)\wedge*_yk(T_zx,y)\,.
\end{align*}
Therefore $C_k$ is translation invariant if $k$ is doubly
translation invariant in the sense that
\[
        k(x,y)=k(S_z x, S_zy)\qquad\forall S_z\,.
\]
Using the translation-invariant basis of $m$\forms\ $w^I$ we see
that the general expression of $k$ is
\[
    k(x,y) =\sum_{I,J} k_{I,J}(x,y) w^I(x)\otimes w^J(y)
\]
with $k_I(x,y)$ doubly-invariant functions, that is, of the form
$k_{I,J}(x,y)=a_{I,J}(S_yx)$ for some function (or distribution)
$a_{I,J}$. If $\delta_0$ denotes the Delta-mass at $e$ and
\[
    \delta(x,y)=\sum_{I,J}\delta_0(S_y x)w^I(x)\otimes w^J(y)
\]
then formally
\[
        \alpha(x)=\int_{\Hn}\alpha(y)\wedge *_y\delta(x,y)\,.
\]
If $P$ is an operator on $m$\forms\ commuting with the $T_y, S_y$,
we will thus have
\[
        P\alpha(x)=\int_{\Hn} \alpha(y)\wedge *_y P_x(\delta(x,y))
\]
and indeed $k(x,y)=P_x(\delta(x,y))$ is formally doubly-invariant.
This shows, in loose terms, that the operator $C_k$ of convolution
with a doubly translation invariant kernel $k$ gives the general
translation-invariant operator acting on $m$\forms. If
$$k(x,y)=\dsum_{I,J} a_{I,J}(S_y x)w^I(x)\otimes w^J(y)$$ and
$\alpha(x)=\sum\alpha_I(x)w^I(x)$, then $C_k \alpha$ has in the
basis $w^I(x)$ coefficients given by
\[
        (C_k\alpha)_I(x)=\sum_J\int_{\Hn} a_{I,J}(S_yx)
        \alpha_J(y)\,d\mu(y)\,.
\]
Thus in the basis $w^I$ everything reduces of course to
convolution of functions. For a function convolution kernel $a(S_y
x)$ and a test function $u\in\dcal(\Hn)$ we may think in
\[
    C_au(x)=\int_{\Hn} u(y) a(S_y x)\,d\mu(y)
\]
as an infinite linear combination of inverse translates $a(S_yx)$
of $a(x)$. Since the vector fields $X_i$ commute with
translations, it follows that whenever everything makes sense
\begin{equation}\label{eq4}
    X_i(C_au) = C_{X_ia} u\,.
\end{equation}
We point out that this convolution is not commutative, $C_au$ is
in general different from $C_u a$. Correspondingly, $X_i C_a
u-C_aX_iu$ is in general not zero; in fact one can easily show
(\cite[lemma 3.1]{BG}) that these commutators are linear
combinations of other convolution operators built from $a(S_yx)$.

\newpage

\noindent {\bf 2.2.} Let $P$ be a generic translation-invariant
operator acting on $m$\forms. We have seen in the previous
subsection that we can associate to $P$ a doubly-translation
invariant kernel $k(x,y)$ so that $P=C_k$. By the same argument as
before, $P$ will be isometry invariant if and only if $k(\varphi
x,\varphi y)= k(x,y)$ $\forall\varphi\in\Iso(\Hn)$ in which case
we say that $k$ is \emph{doubly isometry-invariant}. Working in
the ball model and since every $\varphi\in\Iso(\Hn)$ is the
composition of a translation with some $U\in O(n)$, the additional
requirement on the kernel \linebreak $k(x,y)=\sum
a_{I,J}(S_yx)w^I(x)\otimes w^J(y)$ amounts to $k(Ux,U0)=k(x,0)$,
that is,
\[
    \sum_{I,J}a_{I,J} (Ux)U^*w^I(x)\otimes U^*w^J(0)
        =\sum_{I,J}a_{I,J}(x)w^I(x)\otimes w^J(0), \forall U
\]
Thus we are interested in describing those $k(x,0)$ ---which is a
$m$\form\ at $0$ whose coefficients are $m$-forms in $x$--- that
are doubly invariant by all $U\in O(n)$ in the sense above. Once
the $k(x,0)$ having this property are known, $k(x,y)=k(S_yx,0)$
defines the general doubly isometry invariant $m$\form. For $m=0$
the $k(x,0)$ are simply the radial functions $a(|x|)$, and $a(|S_y
x|) = a(|\varphi_yx|)$ is the general doubly isometry invariant
function. For $m\neq 0$ their general expresion is not so simple.
We find it more convenient to use the usual basis $dx^I$ so we
look at $k(x,0)$ in the form
\begin{equation}\label{eq5}
    k(x,0)= \sum_{|I|=|J|=m} b_{I,J}(x)\,dx^I\otimes dx^J(0)
\end{equation}
and we must impose $\dsum_{I,J} b_{I,J}(Ux) d(Ux)^I\otimes
d(Ux)^J(0)=k(x,0)$ $\forall U$. For instance
$$\gamma(x,0)=\dsum_{i=1}^{n}dx^i\otimes dx^i(0)$$
is easily seen to be doubly $O(n)$\invariant, and so is
$$\gamma_m=\frac{1}{m!}\gamma\wedge\dots\wedge\gamma=\dsum_{|I|=m}
dx^I\otimes dx^I(0)$$ (here we use the symbol $\wedge$ to denote
as well the exterior product of double forms defined by
$(\alpha_1\otimes\beta_1)\wedge(\alpha_2\otimes\beta_2)=
(\alpha_1\wedge\alpha_2)\otimes(\beta_1\wedge\beta_2)$). Another
doubly $O(n)$\invariant\ 1\form\ is
$$\tau(x,0)=\left(\dsum_{i=1}^n
x_i\,dx^i\right)\otimes \left(\dsum_{i=1}^n x_i\,dx^i(0)\right).$$
\begin{lemma}\label{lem1}
The double forms $\gamma$ and $\tau$ generate all doubly
$O(n)$\invariant\ $k(x,0)$. More precisely, their general
expression in the ball model is
\begin{align}\label{eq6}
  k(x,0)&= A_1(|x|)\gamma_m + A_2(|x|)\tau\wedge\gamma_{m-1}\,
  \qquad 0<m<n\\[2mm]
  k(x,0)&= A(|x|) \gamma_m\,,\qquad m=0,n\nonumber
\end{align}
\end{lemma}

\begin{proof}
First we prove by induction the following statement $S(n)\colon$
if $k(x,0)$ is a doubly invariant $(p,q)$\form\
$\dsum_{|I|=p,|J|=q} c_{I,J}\,dx^I\otimes dx^J(0)$ with constant
coefficients, then $k\equiv 0$ if $p\neq q$, or $k$ is diagonal
i.e. $k(x,0)=c\dsum_{|I|=p} dx^I\otimes dx^I(0) = c\gamma_p$ if
$p=q$.
\noindent Of course $S(1)$ is obvious; assuming $S(n-1)$,
let us break $k(x,0)$ in four pieces, depending on whether
$i_1,j_1=1$ or not:
\[
    k=\sum_{i_1=j_1=1} c_{I,J}\,dx^I\otimes dx^J(0)
    +\sum_{i_1=1, j_1\neq1}
    +\sum_{i_1\neq 1, j_1=1}
    +\sum_{i_1\neq 1, j_1\neq 1}\stackrel{\text{def}}{=} k_1+k_2+k_3+k_4\,.
\]
We may write $k_1=(dx^1\otimes dx^1(0))\wedge \widetilde{k_1}$,
$k_2=(dx^1\otimes 1)\wedge \widetilde{k_2}$, $k_3=(1\otimes
dx^1(0))\wedge \widetilde k_3$, with $\widetilde{k_1}$,
$\widetilde{k_2}$, $\widetilde{k_3}$, $k_4$  double forms in the
$dx^2,\dots, dx^n, dx^2(0),\dots, dx^n(0)$ of bidegrees $(p-1,
q-1)$, $(p-1, q)$, $(p,q-1)$ and $(p,q)$, respectively. Imposing
that $k$ is doubly invariant by $U$ of the type
\begin{equation}\label{eq7}
    U=\begin{pmatrix}
       1&0&0&\cdots\\
       0&&&\\
       0&& U_1&\\
       \vdots &&&\\
       0 &&&
       \end{pmatrix}\,,\qquad U_1\in O(n-1)
\end{equation}
we see that $\widetilde{k_1}$, $\widetilde{k_2}$,
$\widetilde{k_3}$ and $k_4$ are $O(n-1)$\invariant. We apply the
induction hypothesis: if $p=q$
$\widetilde{k_2}=\widetilde{k_3}=0$, and $\widetilde{k_1}, k_4$
are diagonal, i.e.
\begin{equation}\label{eq8}
    k=c_1 \sum_{i_1=1} dx^I\otimes dx^I(0) + c_2\sum_{i_1\neq 1}
    dx^I\otimes dx^I(0)\,.
\end{equation}
If we use now $U\in O(n)$ permuting the first two axes we see that
$c_1=c_2$ and hence $k$ is diagonal, establishing $S(n)$ in case
$p=q$. If $|p-q|>1$ everything is $0$. Finally if $|p-q|=1$, say
$p=q+1$, then $\widetilde{k_2}$ is diagonal and all others are
zero
\[
  k=c(dx^1 \otimes 1)\wedge \sum_{|J|=q} d{x'}^J \otimes d{x'}^J(0)
\]
where $x'=(x_2,\dots x_n)$. If we impose the invariance under the
permutation of the first two axes as before, it is clear that $k$
must be zero.

Having proved that $S(n)$ holds for all $n$, let now $k(x,0)$ be
as in \eqref{eq5} doubly $O(n)$\invariant. Clearly $k(x,0)$ is
then determined by its values $k(\rvec, 0)$, where
$\rvec=(r,0,0,\dots,0)$. Fixed $r$, $k(\rvec,0)$ may be regarded
as a double $(m,m)$\form\ with constant coefficients, which is
invariant by all $U\in O(n)$ fixing $\rvec$, that is, of type
\eqref{eq7}. We write now the decomposition of $k(\rvec,0)$ in
terms of $\widetilde{k_1}(r,0)$,  $\widetilde{k_2}(r,0)$,
$\widetilde{k_3}(r,0)$ and $k_4(r,0)$ as before, and applying
$S(n)$ we get \eqref{eq8}
\[
        k(\rvec,0)= c_1(r)\sum_{i_1=1} dx^I(\rvec)\otimes dx^I(0)
        +c_2(r)\sum_{i_1\neq 1} dx^I(\rvec) \otimes dx^I(0)
\]
(if $m=n$ the last term is zero and the first is $\gamma_m$) which
we write
\begin{align*}
  &= (c_1(r)-c_2(r))\sum_{\substack{i_1=1\\ |I|=m}} dx^I(\rvec)
  \otimes dx^I(0) + c_2(r)\sum_{|I|=m} dx^I(\rvec)\otimes
  dx^I(0)\\[2mm]
  &= (c_1(r)-c_2(r)) dx^1(\rvec)\otimes dx^1(0)\wedge
  \sum_{|I|=m-1} dx^I(\rvec)
  \otimes dx^I(0) + c_2(r)\sum_{|I|=m} dx^I(\rvec)\otimes
  dx^I(0)\\[2mm]
  &= (c_1(r)-c_2(r))r^{-2} \tau(\rvec,0) \gamma_{m-1}(\rvec, 0) +
  c_2(r) \gamma_m(\rvec, 0)\,.
\end{align*}
Finally, fixed $x$ we choose $U$ such that $Ux=\rvec$, $r=|x|$,
and use the invariance of $k$, $\tau$, $\gamma$ to find
\eqref{eq6}with $A_1(r)=c_2(r)$, $A_2(r)=r^{-2}(c_1(r)-c_2(r))$.
\end{proof}

\vspace{1cm} To find the general expression of a doubly isometry
invariant kernel $k(x,y)$ we must translate $k(x,0)$ to an
arbitrary point: $k(x,y)= k(S_y x, S_yy)$. We may use any isometry
mapping $y$ to $0$, for instance we may use $\varphi_y$ given by
\eqref{eq1} instead of $S_y$. We introduce the basic forms
$\alpha,$ $\beta,$ $\tau$ and $\gamma$
\begin{align*}
  \alpha &= \alpha(x,y) =\sum_{i}\varphi_y^i(x)\,d\varphi_y^i(x)\,,
  \qquad \beta=\sum_i\varphi_y^i(x)\,d\varphi_y^i(y)
  =-\sum_i\varphi_y^i(x) \frac{dy^i}{1-|y|^2}\\[3mm]
  \tau&= \alpha\otimes\beta;
  \gamma(x,y)= \sum_{i=1}^n d\varphi_y^i(x)\otimes
  d\varphi_y^i(y)= \frac{-1}{1-|y|^2}\sum_{i=1}^n
  d\varphi_y^i(x)\otimes dy^i =d_x\beta\,.
\end{align*}
The lemma gives part (a) of the following theorem. Part (b) gives
other equivalent general expressions, which are intrinsic, that
is, independent of the model of $\Hn$ at use.

\begin{theorem}\label{eqth1} (a) The general expression of an $(m,m)$\form\ $k(x,y)$ doubly
isometry\invariant\ in $\Hn$, in the ball model, is
\begin{align*}
  k(x,y)&= A_1(|\varphi_yx|)\gamma_m(x,y) + A_2(|\varphi_y
  x|)\,\tau(x,y)\wedge\gamma_{m-1}(x,y)\,,\qquad 0<m<n\\[2mm]
  k(x,y)&= A(|\varphi_yx|) \gamma_m(x,y)\qquad m=0,n\,.
\end{align*}
(b) Another equivalent expression for $0<m<n$ is
\begin{align*}
  k(x,y)&= B_1(D) (d_x d_y D)^m + B_2(D)(d_xD\otimes d_yD)
         \wedge(d_xd_yD)^{m-1} =\\[2mm]
        &= (C_1(D)d_xd_y D + C_2(D)d_x D\otimes d_yD)^m
\end{align*}
where $D$ denotes an arbitrary function of the geodesic distance
$d(x,y)$.

\noindent (c) All such $k(x,y)$ are symmetric in $x,y\in \Hn$
\end{theorem}

\begin{proof}
Part (a) has been already proved. For (b) note first that it is
enough to consider \emph{one} function of $d$: we choose
$D=r(x,y)^2$ which in the ball model equals $|\varphi_y(x)|^2$.
Then $d_xD=2\alpha$, and using \eqref{eq1}, \eqref{eq2} one finds
\[
        d_y D = 2(1-D) \sum_i\varphi_y^i(x) \frac{d_y^i}{1-|y|^2}
        =-2(1-D)\beta\,.
\]
This gives $\tau=\alpha\otimes\beta
=\frac{-1}{4}\frac{1}{1-D}d_xD\otimes d_y D$, and
\[
        d_x d_y D = +2 d_xD\otimes\beta -2(1-D)d_x \beta =
        +4\tau-2(1-D)\gamma\,.
\]
Therefore $(d_xd_yD)^{m-1}$ and $2^{m-1}(1-D)^{m-1}\gamma_{m-1}$
differ in a term containing $\tau$, and so (b) follows. Part (c)
is a consequence of (b).
\end{proof}

We will need the expression of the generators $\tau$, $\gamma$ in
terms of the invariant basis $w^i$. We obtain these using formula
\eqref{eq2} for $r^2(x,y)$ in the half-space model. First
\[
\alpha=\frac{d_x r^2}{2}=\frac{1-r^2}{2(|x-y|^2+
    4x_ny_n)}\left(2\sum_{i=1}^{n-1} x_n(x_i-y_i)w^i(x)+
    (2x_n(x_n-y_n)-|x-y|^2)w^n(x)\right)
\]
\[
\beta=\frac{d_yr^2}{2(r^2-1)}=\frac{-1}{2(|x-y|^2+4x_ny_n)}\left(2\sum_{j=1}^n
y_n(y_j-x_j)w^j(y)+(2y_n(y_n-x_n)-|x-y|^2)w^n(y)\right).
\]
In the following we write $w^{ij}=w^i(x)\otimes w^j(y)$. We have
\[
        \tau=\alpha\otimes\beta=\frac 14 \frac{1-r^2}{(|x-y|^2 + 4x_ny_n)^2}
        \sum_{ij}P_{i,j}(x,y)w^{i,j}
\]
where the $P_{ij}(x,y)$ are certain homogeneous polynomials. As we
know, everything can be written in terms of $z=S_yx$: for instance
\[
        1-r^2 = \frac{4x_ny_n}{|x-y|^2 + 4x_ny_n} = \frac{4z_n}{|z|^2+2z_n + 1}
\]
and say for $i,j<n$
\[
        \frac{P_{ij}}{(|x-y|^2+ 4x_ny_n)}=
        \frac{x_ny_n(x_i-y_i)(x_j-y_j)}{(|x-y|^2 + 4x_ny_n)^2} =
        \frac{z_n z_iz_j}{(|z|^2+2z_n+1)^2}\,.
\]
Therefore we may write
\begin{equation}\label{eq9}
    \tau=\frac{1-r^2}{(|z|^2+2z_n+1)^2}\sum_{i,j} p_{i,j}(z) w^{i,j}\,.
\end{equation}
For $\gamma=d_x\beta$ we obtain a similar expression
\begin{multline*}
\frac{4}{1-r^2}\,\gamma=\sum_{i,j=1}^{n-1}\left(\delta_{ij}-\frac{2(x_i-y_i)(x_j-y_j)}{|x-y|^2
+ 4x_ny_n}\right)w^{i,j}+\\
+\left(1-\frac{2\dsum_{i=1}^{n-1}|x_i-y_i|^2}{|x-y|^2+4x_ny_n}\right)
w^{n,n}+\sum_{i=1}^{n-1}\frac{2(x_i-y_i)(x_n-y_n)}{|x-y|^2 +
4x_ny_n}(w^{i,n}-w^{n,i}).
\end{multline*}
Again this can be written
\begin{equation}
    \gamma=\frac{1-r^2}{4(|x-y|^2 + 4x_ny_n)} \sum_{i,j}Q_{ij}(x,y)
        w^{i,j}=\frac{1-r^2}{(|z|^2+2z_n +1)}
            \sum q_{ij}(z)w^{i,j}\,.
\end{equation}
Notice that
\[
            \frac{p_{ij}(z)}{(|z|^2 + 2z_n + 1)^2}= O(1)\,,
      \qquad \frac{q_{ij}(z)}{(|z|^2 + 2z_n+1)} =O(1)
\]
and hence
\begin{equation}\label{eq11}
        |\tau(x,y)|=O(1-r^2)\,,\qquad |\gamma(x,y)|=O(1-r^2)\,.
\end{equation}

\section{Riesz forms and Riesz form-potentials in $\Hn$}

\noindent {\bf 3.1.} Our next objective is now to find an explicit
left-inverse $L$ for $\Delta$ on $\dcal_m(\Hn)$. Since $\Delta$ is
invariant by all isometries, $L$ should be too. By what has been
discussed in section 2, $L$ should have a kernel $k_m(x,y)$,
\[
    L\eta(x)=\int_{\Hn}\eta(y)\wedge *_y k_m(x,y)
\]
doubly invariant by all isometries. Alternatively, notice that if
$k$ is \emph{some} kernel such that
\begin{equation}\label{eq12}
    \eta(x)=\int_{\Hn} \Delta\eta(y)\wedge
    *_yk(x,y)\,,\qquad\eta\in\dcal_m(\Hn)
\end{equation}
(which formally exists because $\Delta\eta =0$,
$\eta\in\dcal_m(\Hn)$ imply $\eta=0$) then its average over the
unitary group $O(n)$ with respect the normalized left-invariant
measure $d\mu(U)$,
\[
        k_1(x,y)=\int_{O(n)} k_0(Ux,Uy)d\mu(U)
\]
still satisfies \eqref{eq12} and it is doubly invariant by $O(n)$.
If $\varphi_x$ is an isometry mapping $x$ to $0$,
$k_2(x,y)=k_1(\varphi_x x,\varphi_x y)$ is independent of
$\varphi_x$, satisfies \eqref{eq12} and is doubly invariant by all
isometries.

Anyway, we look for a doubly isometry-invariant kernel $k_m$ for
which \eqref{eq12} holds, and then consider the operator $L$
defined by $k_m$ as above. Taking for granted by now that this
operator $L$ is well defined on $\dcal_m(\Hn)$ and maps
$\dcal_m(\Hn)$ into locally integrable $m$\forms, notice that
\eqref{eq12} and the symmetry of $k_m$ together imply that $L$ is
a right-inverse too, that is, $\Delta L\alpha=\alpha$ for
$\alpha\in\dcal_m(\Hn)$ in the weak sense:
\begin{align*}
     \langle \Delta L\alpha,\eta\rangle
    & = \langle L\alpha,
        \Delta\eta\rangle =\int_x L\alpha(x)\wedge \aster  \Delta\eta(x)=\\[2mm]
    &=\int_x
        \left\{\int_y \alpha(y)\wedge *_yk_m(x,y)\right\}
    \wedge \aster \Delta\eta(x)\\[2mm]
    & =\int_y \alpha(y)\wedge *_y
        \left\{\int_x k_m(x,y)\wedge \aster \Delta\eta(x)\right\}
        =\langle \alpha,\eta\rangle\,.
\end{align*}
We work in the ball model. By Theorem 2.2, $k_m(x,y)$ is of type
\begin{align*}
        k_m(x,y)&= A(|\varphi_xy|)\gamma_m\,, \quad m=0,n\\[3mm]
        k_m(x,y)&= A_1(|\varphi_x y)\gamma_m + A_2(|\varphi_xy|)
        \tau\wedge \gamma_{m-1}\,,\qquad 0<m<n
\end{align*}
where $\gamma\!=\!\dsum_i d\varphi_x^i(x)\otimes d\varphi_x^i(y),
\tau\!=\! \alpha\otimes\beta$ with
$\alpha\!=\!\dsum_i\varphi_x^i(y) d\varphi_x^i(x)$,
$\beta\!=\!\dsum_i\varphi_x^i(y) d\varphi_x^i(y)$ (notice that we
are exchanging $x,y$, using (c) in Theorem 2.2). Condition
\eqref{eq12} implies $\Delta_y k_m(x,y) =0$ in $y\neq x$ (while
$\Delta Lw=w$ implies $\Delta_xk_m(x,y)=0$ in $x\neq y$). In fact,
\eqref{eq12} amounts to requiring $\Delta_yk_m(x,y)=\delta_x$ in a
sense to be described below.

\vspace{1cm}

\noindent {\bf 3.2.} In a first step we look for conditions on the
$A_1,$ $A_2,$ so that $\Delta_yk_m(x,y)=0$ in $y\neq x$. A lengthy
computation will show that the general harmonic $k_m$ depends on
four parameters. By the invariance of $k_m$, we may assume $x=0$,
in which case, writing $r=|y|$,
\begin{align*}
        k_m(x,y)&= A(r)\gamma_m\,,\qquad m=0,n\\[2mm]
        k_m(0,y)&= A_1(r)\gamma_m + A_2(r)\tau\wedge \gamma_{m-1}
\end{align*}
with $\gamma=\sum dx^i(0)\otimes dy^i$, $\tau=\alpha\otimes\beta,
\alpha=\sum y^i dx^i(0)$, $\beta=rdr$. Since $*_x*_yk_m(x,y)$ is
again doubly invariant, it must have an analogous expression with
$m$ replaced by $n-m$. Indeed, it is easily checked that
\begin{align*}
        *_x *_y\gamma_m
    &= \frac{m!}{(n-m)!}
        (1-r^2)^{2m-n}\gamma_{n-m}\\[3mm]
        *_x*_y (\tau\wedge \gamma_{m-1})
    &= (m-1)!(1-r^2)^{2m-n}
        \left(r^2\frac{\gamma_{n-m}}{(n-m)!} -
            \frac{\tau\wedge\gamma_{n-m-1}}{(n-m-1)!}\right)
\end{align*}
whence
$*_x*_yk_m(x,y)=\frac{m!}{(n-m)!}(1-r^2)^{2m-n}\gamma_{n-m}$ for
$m=0,n$ and for $0<m<n$,

\begin{equation}\label{eq13}
        *_x*_yk_m(0,y) =\\[2mm]
        \frac{(m-1)!(1-r^2)^{2m-n}}{(n-m)!}
        \left[(mA_1+ r^2A_2)\gamma_{n-m}
        -(n-m)A_2\tau\wedge\gamma_{n-m-1}\right].
\end{equation}
Moreover, since $*$ commutes with $\Delta$, it is natural to
require as well that $*_x*_y k_m=k_{n-m}$, that is, we may assume
from now on that $0\leq m\leq n/2$.

\noindent For $m=0$, using \eqref{eq3.3} we find
\[
        \Delta(A(r))=\frac14 (1-r^2)\left[-1(1-r^2)A''+((3-n)r + r^{-1}(1-n))A'\right]
\]
from which it follows that $A'(r)=c_0(1-r^2)^{n-2} r^{1-n}$ and
\[
        A(r)=c_1-c_0\int_r^1(1-s^2)^{n-2}s^{1-n}\,ds\,.
\]
We start now computing $\Delta_y k_m(0,y)$ for $0<m \leq n/2$
using that on $m$\forms $\,\Delta$ equals $(-1)^{m+1}(\aster
d\aster d+(-1)^n d\aster d\aster )$. The double form $\Delta_y
k_m(x,y)$ is also doubly invariant, and therefore it must have the
same expression as $k_m$ with $A_1$, $A_2$ replaced by other
functions $B_1$, $B_2$ to be found. In the computations we will
use besides \eqref{eq13} the equations
\begin{align*}
       & d_y\alpha =\gamma\,,\quad d_y(\tau\wedge \gamma_{m-1})
         =-rdr\wedge\gamma_m = -\beta \wedge \gamma_m\\[3mm]
       & *_x*_y dr\wedge \gamma_m =(-1)^m\frac{m!}{(n-m-1)!} (1-r^2)^{2m+2-n}
        r^{-1}\alpha\wedge \gamma_{n-m-1}
\end{align*}
which are easily checked as well. First, $d_y
k_m(0,y)=(A'_1-rA_2)dr\wedge \gamma_m$, so by the equations above
\begin{align}\label{eq14}
  *_x*_yd_yk_m(0,y)=&(-1)^m\frac{m!}{(n-m-1)|!} (1-r^2)^{2m+2-n}
        (A'_1-rA_2)r^{-1}\alpha\wedge \gamma_{n-m-1}=\\[2mm]
    &\stackrel{\text{def}}{=}\frac{(-1)^m
    m!}{(n-m-1)!}A_3\alpha\wedge \gamma_{n-m-1}\,,\nonumber
\end{align}
\[
        *_xd_y*_yd_yk_m(0,y)=\frac{(-1)^mm!}{(n-m-1)!}
        \left(A_3\gamma_{n-m} +A'_3r^{-1}\tau \wedge\gamma_{n-m-1}\right)
\]
\begin{align*}
& *_yd_y*_yd_yk_m(0,y)=(-1)^{m(n-m-1)}*_y*_x(A_3\gamma_{n-m}+
    A'_3r^{-1}\tau\wedge \gamma_{n-m-1})=\\[2mm]
& =(-1)^{m(n-m-1)}\frac{m!}{(n-m-1)!}(1-r^2)^{n-2n}
    \left( A_3\frac{(n-m)!}{m!}\gamma_m+
            A'_3 r\frac{(n-m-1)!}{m!}\gamma_m\right. \\[2mm]
    &\qquad\qquad \left.- A'_3 r^{-1}\frac{(n-m-1)!}{(m-1)!}
                \tau\wedge\gamma_{m-1}\right)\\[2mm]
&=(-1)^{m(n-m+1)}(1-r^2)^{n-2m} \left[((n-m)A_3 + A'_3r)\gamma_m -
mA'_3r^{-1}\tau\wedge \gamma_{m-1}\right]\,.
\end{align*}
By the analogous computation, applying $d_y$ to \eqref{eq13}
\begin{align*}
*_xd_y*_yk_m(0,y)&=\frac{(m-1)!}{(n-m)!}\left[\left[(mA_1+r^2
A_2)(1-r^2)^{n-2m}\right]'\right.\\[2mm]
& \qquad\left.+ (n-m)rA_2(1-r^2)^{n-2m}\right] dr\wedge
\gamma_{n-m}
\end{align*}
\begin{align}\label{eq15}
 &*_yd_y*_yk_m(0,y)=(-1)^{(m+1)(n-m)}(1-r^2)^{2m+2-n}r^{-1}\nonumber\\[2mm]
 &\qquad=\left[\left[ (mA_1+r^2A_2)(1-r^2)^{n-2m}\right]'+
    (n-m)rA_2(1-r^2)^{n-2m}\right]\alpha\wedge\gamma_{m-1}=\\[2mm]
 & \qquad\stackrel{\text{def}}{=}(-1)^{(m+1)(n-m)}  A_4\alpha\wedge \gamma_{m-1}\nonumber
\end{align}
\[
        d_y*_yd_y*_yk_m(0,y)=(-1)^{(n-m)(m+1)}(A'_4r^{-1}\tau\wedge\gamma_{m-1}+ A_4\gamma_m)
\]
It follows finally that $\Delta=(-1)^{nm+1}(\aster d\aster
d+(-1)^n d\aster d\aster )$ on $k_m$ equals
\[
    \Delta_yk_m(0,y)= B_1\gamma_m + B_2\tau\wedge\gamma_{m-1}
\]
with
\begin{align*}
    B_1&=- A_4 -(1-r^2)^{n-2m}((n-m)A_3+A'_3r)\\[2mm]
    B_2&=-A_4 r^{-1}+ m(1-r^2)^{n-2m}A'_3r^{-1}\,.
\end{align*}
Therefore, $\Delta_yk(0,y)=0$ is equivalent to the system $B_1=0$,
$B_2=0$. It easily follows from this that $A_3$ satisfies the
equation
\[
        r(1-r^2)A''_3 + \left[(n+1)-r^2(3n+1-4m)\right]A'_3
        -2(n-2m)(n-m)rA_3=0\,.
\]
Replacing in the equation $B_1=0$, $A_4$ by its expression in
terms of $A_1$ and $A_2$, and then $A_2$ by its expression in
terms of $A_1$ and $A_3$, we find that $A_1$ satisfies the
inhomogeneous equation
\begin{multline*}
        r(1-r^2)A''_1+ \left[(n+1)+ (n-1-4m)r^2\right]A'_1
            +2m(n-2m)rA_1 = \\[2mm]
            = 2rA_3(1+r^2)(1-r^2)^{n-2m-2}
\end{multline*}
The change of variables $A_1(r)=G(x)$, $A_3(r)=H(x)$, $x=r^2$,
transforms these into the hypergeometric equations
\begin{align}\label{eq16}
    x(1-x)H''(x)&+\left[\frac n2 + 1 -(\frac 32n+1-2m)x\right]
        H'(x)-(\frac n2-m)(n-m)H=0\nonumber\\[2mm]
   x(1-x)G''(x)&+\left[\frac{n}{2}+1-(2m+1-\frac
   n2)x\right]G'(x)-m(m-\frac n2)G=\\[2mm]
   &\qquad\qquad=\frac
   12(1+x)(1-x)^{n-2m-2}H(x)\stackrel{\text{def}}{=}f(x)\nonumber
\end{align}
This system is equivalent to $\Delta_y k_m(x,y)=0$ in $y\neq x$,
whence the general doubly-invariant $k_m$ harmonic in $y\neq x$
depends on four parameters. Note that for $m=\frac n2$ the
homogeneous equations are the same and can be solved explicitely:
the general solution is $H=as^{-\frac n2}+ b$ and
\begin{align}\label{eq17}
    G(x)= cx^{-\frac n2} + d + \frac 12\int_{1/2}^x t^{-\frac n2-1}
    \left\{ \int_{0}^t s^{n/2}(1+s)(1-s)^{-3}(as^{-\frac n2} +
    b)\,ds\right\}\,dt\,.
\end{align}
For $m<\frac n2$ a fundamental family for the first equation in
(16) is given by
$$u_1(x)= x^{-\frac n2} F(-m, \frac n2-m, 1-\frac n2, x),
u_2(x)= F(\frac n2 -m, n-m, \frac n2+1, x)$$ The hypergeometric
function in $u_1$ is a polynomial in $x$ of degree $m$ with
positive coefficients, $1+x$ if $m=1$. A fundamental family for
the second is given by
\begin{multline*}
u_3(x)= x^{-\frac n2}F(m-n,m-\frac n2, 1-\frac n2, x)= x^{-\frac
n2}(1-x)^{n+1-2m} F(\frac n2+1-m, 1-m, 1-\frac n2, x),\\
 u_4(x) =
F(m, m-\frac n2, 1+\frac n2,
x)\quad\quad\quad\quad\quad\quad\quad\quad\quad\quad\quad\quad
\end{multline*}
The  hypergeometric function in $u_3$ is a polynomial of degree
$m-1$ with positive coefficients (see \cite{E} for all these
facts). The wronskian $w(x)$ for this second equation is, by
Liouville's formula
\[
    W(x)=W(x_0)\exp-\int_{x_0}^x\frac{\frac n2 +1-(2m+1-\frac n2)t}
                {t(1-t)}\,dt = c_{mn} x^{-\frac n2-1}(1-x^{n-2m})
\]
It follows from this that the parametrization for $G$ is given by
\begin{equation}\label{eq18}
        G(x)= c(x)u_3(x)+ d(x) u_4(x)
\end{equation}
where $c(x), d(x)$ satisfy, with $H(x)= a u_1(x)+bu_2(x)$,
\begin{align*}
        c'(x)&= \frac{u_4(x)f(x)}{x(1-x)W(x)}
            =\frac 12 \, c_{mn}^{-1} H(x)(1+x)x^{\frac n2}(1-x)^{-3}u_4(x)\\
                    d'(x)&= -\frac{u_3(x)f(x)}{x(1-x)W(x)}
            =-\frac 12\, c_{mn}^{-1} H(x)(1+x)x^{\frac n2}(1-x)^{-3}
            u_3(x)\,.
\end{align*}
Once $A_1(r)= G(r^2)$ and $A_3(r)=H(r^2)$ are known, the kernel
$k_m(x,y)$ is completely known, because by the definition of $A_3$
in \eqref{eq14}, $A_2(r)=-(1-r^2)^{n-2m-2} A_3(r)+ r^{-1}A'_1(r)=
-(1-x)^{n-2m-2}H(x)+ 2G'(x)$.

The choice $a=0$, $c(0)=0$ ($a=c=0$ in the parametrization
\eqref{eq17} for $m=\frac n2$) gives all doubly invariant
$k_m(x,y)$ which are \emph{globally} harmonic, with no
singularity, and they are therefore spanned by the forms
corresponding to the choice $G=u_4$ and to the choice $a=0$,
$b=1$, $c(0)=0$, $d(0)=0$,
\begin{multline*}
        G(x)=\left\{\int_0^x(1+t)(1-t)^{-3} t^{n/2} u_2(t)u_4(t)\,dt\right\}
        u_3(x)\\[3mm]
        -\left\{\int_0^x (1+t)(1-t)^{-3}
        t^{n/2}u_2(t)u_3(t)\,dt\right\}u_4(x)\,.
\end{multline*}
As a particular case, note that for $m=\frac n2, \gamma_m$ is
harmonic in $\Hbb^{2m}$, and it is the simplest example of a
non-zero harmonic $m$-form in $L^2(\Hbb^{2m})$.

\vspace{1cm}

\noindent {\bf 3.3.} Besides being harmonic in $y\neq x$, the
singularity at $y=x$ most be such that \eqref{eq12} holds. Again,
we may assume $x=0$; we check this property using \emph{second's
Green identity}, whose version for general forms we recall now.

The operator $\delta$ being the adjoint of $d$, one has for a
smooth domain $\overline{\,\Omega}\subset \Bn$ and $\alpha, \beta$
smooth forms on $\overline{\,\Omega}$ with $\text{degree} \alpha=
\text{degree}\beta-1$
\[
        \int_{\partial\Omega}\alpha\wedge \aster \beta =\int_{\Omega}
        d\alpha\wedge\aster \beta -\int_\Omega
        \alpha\wedge\aster \delta\beta\,.
\]
Given two $m$\forms\ $\eta$, $\omega$, applying this with
$\alpha=\delta\eta$, $\beta=\omega$, next with $\alpha=\omega$,
$\beta=d\eta$ and substracting one gets \emph{the first Green's
identity for m-forms}
\[
    \int_{\partial\Omega}(\delta\eta\wedge\aster \omega-\omega\wedge \aster d\eta)
    =\int_\Omega(\Delta\eta\wedge \aster \omega-\delta\eta\wedge \aster \delta\omega-d\eta\wedge \aster d\omega)
\]
Permuting $\omega, \eta$ and substracting again gives \emph{the
second Green's identity}
\[
    \int_{\partial\Omega}(\delta\eta\wedge \aster  \omega-\omega \wedge \aster d\eta
    -\delta\omega\wedge\aster \eta + \eta\wedge \aster  d\omega)
    =\int_\Omega(\Delta\eta\wedge\aster \omega-\Delta\omega\wedge
    \aster \eta)\,.
\]
We apply this to $\Omega=B(0,R)-B(0,\varepsilon)$
$0<\varepsilon<R<1$, $\eta\in\dcal_m(\Hn)$ and our $k_m(0,y)$ to
get
\begin{equation}\label{eq19}
    \int_{|y|\geq\varepsilon}\Delta\eta\wedge *_y k_m(0,y)
    =\int_{|y|=\varepsilon}(k_m\wedge \aster  d\eta + \delta_yk_m\wedge \aster \eta
        -\delta\eta\wedge \aster _yk_m - \eta\wedge \aster  dk_m)\,.
\end{equation}
In case $m=0$, the terms in $\delta k_m, \delta\eta$ are of course
zero; to get a term in $\eta(0)$ on the right when $\varepsilon\to
0$ we need $dk_m$ of the order of $\varepsilon^{1-n}$ and $k_m$ of
the order of $\varepsilon^{2-n}$ in $|y|=\varepsilon$. That makes
$k_m$ locally integrable too, and \eqref{eq12} is obtained leting
$\varepsilon\to 0$. This means that for $m=0$ $k$ is unique and is
given by the well-known Green's function
\begin{equation}\label{eq20}
    A(r)=c_n\int_r^1 (1-s^2)^{n-2}s^{1-n}\, ds
\end{equation}
for an appropriate choice of $c_n$. In case $m>0$, again we need
$|k_m(0,y)|=o(r^{1-n})$ as $r\to 0$, so that the first and third
terms on the right have limit $0$ as $\varepsilon\to 0$; then
$k_m$ is integrable in $y$ and the integral on the left converges
to $\int \Delta\eta\wedge\aster k_m$. Using the expression for
$\aster dk_m$ in \eqref{eq14}, we find
\[
    \int_{|y|=\varepsilon}\eta\wedge \aster d_yk_m =
        \frac{(-1)^{m(n-m+1)}n!}{(n-m-1)!}\, A_3(\varepsilon)*_x
        \int_{|y|=\varepsilon}\eta\wedge\alpha\wedge\gamma_{n-m-1}\,.
\]
By Stoke's theorem, and since $\alpha=O(r)$, the last integral
equals
\[
        (-1)^m\int_{|y|<\varepsilon} \eta\wedge\gamma_{n-m} +
        O(\varepsilon)\,.
\]
If $A_3(\varepsilon)= a_0\varepsilon^{-n}+\dotsc$, we see that
\[
        \lim_\varepsilon\int_{|y|=\varepsilon}\eta\wedge\aster d_y k_m
        = c_n(n-m)\, m! \,a_0 \eta(0)\,.
\]
Using \eqref{eq15} for $\delta k_m=(-1)^{n(m+1)+1}\aster d\aster
$, and proceeding in the same way,
\begin{align*}
    \int_{|y|=\varepsilon} \delta_yk_m\wedge \aster  \eta
    & = -A_4 (\varepsilon)\int_{|y|=\varepsilon} \alpha\wedge
    \gamma_{m-1}\wedge \aster \eta\\[2mm]
    & = -A_4(\varepsilon) \int_{|y|<\varepsilon}(\gamma_m\wedge \aster \eta +
    O(\varepsilon))\,.
\end{align*}
But by the equation $B_1=0$,
$A_4(\varepsilon)=-(1-\varepsilon^2)^{n-2m}
((n-m)A_3(\varepsilon)+\varepsilon A'_3(\varepsilon))=
a_0m\varepsilon^{-n} + O(\varepsilon^{1-n})$, and hence the limit
of the above expression is $-c_n m!a_0 m\eta(0)$. Altogether, we
conclude that if $A_3(\varepsilon)= a_0\varepsilon^{-n} +
o(\varepsilon^{1-n})$ and $k_m(0,y)=o(r^{1-n})$, one has
\[
    \int\Delta\eta\wedge *_yk_m(0,y) = -c_n\, n\, m!\, a_0\, \eta(0)
\]
so \eqref{eq12} will hold for an appropriate choice of $a_0$.
Taking into account the definition of $A_3$ in \eqref{eq14} and
that $|k_m|\simeq |A_1| + r^2|A_2|$, we see from \eqref{eq17} that
if $m=\frac n2$ this is accomplished by the choice $c=0$, $a=a_0$;
then $G(x)\sim\log x$, $A_1(r)\sim \log r,\, A'_1(r)=O(1/r)$,
$A_2(r)=O(r^2)$ if $n=2$; if $n>2$ ,$A_1(r)\sim r^{2-n}$ and
$A_2=O(r^{-n})$. For $0<m<\frac n2$, in terms of the functions $H,
G$ introduced before, this translates to $H(x)\sim a_0 x^{-\frac
n2}, G(x)\sim x^{1-\frac n2}$. Now look at the general expression
of $H, G$ in \eqref{eq18}. The condition $H(x)\sim c_0x^{-\frac
n2}$ fixes $a=a_0$; then near $x=0$ $c'(x)$ is bounded and $d'(x)$
behaves like $x^{-\frac n2}$. Since $u_4(x)$ is bounded, the term
$d(x)u_4(x)$ behaves like $x^{1-\frac n2}$. So, we must normalize
$c(x)$ by $c(0)=0$, so that $c(x)=O(x)$ and the other term
$c(x)u_3(x)$ will behave like $x^{1-\frac n2}$.

In conclusion, all this discussion shows that the doubly invariant
kernels $k_m(x,y)$ satisfying \eqref{eq12} constitute a \emph{two
parameter family} described by $H= a_0 u_1(x)+ bu_2(x)$, $c(0)=0$.
The two parameters are $b$ and the constant of integration for
$d(x)$ in \eqref{eq18}. Equivalently, they are obtained by adding
to the form corresponding to $H=a_0u_1(x)$, $c(0)=0$ and say
$d(\frac 12)=0$ the general globally smooth one described before.

\vspace{1cm}

\noindent {\bf 3.4.} In order to produce the best estimates, in a
sense we need to choose the best of the kernels $k_m$. Naturally
enough, we choose the $k_m$ having the best behaviour at infinity,
$x=1$, that is, so that $G, H$ have the best decrease in size as
$x\to 1$. In case $m=\frac n2$, where we already have the
normalization $c=0$, $a=a_0$, the choice $b=-a$ gives the best
growth $H(x)=O(1-x)$ and $G(x)=O(\log(1-x))$.

The hypergeometric function $u_3$ behaves like $(1-x)^{n+1-2m}$
near $x=1$ while $u_4(x)=F(m,m-\frac n2, 1+\frac n2, x)$ is
bounded because $1+\frac n2-m-(m-\frac n2)= 1+n-2m>0$. Similarly,
$u_1$ is bounded near $x=1$; for $u_2(x)=F(\frac n2-m, n-m, \frac
n2+1,x)$ we have $\frac n2+1-(\frac n2-m)-(n-m)=2m+1-n$ and hence
it behaves like $(1-x)^{2m+1-n}$ if $2m<n-1$ and like $\log(1-x)$
if $2m=n-1$. We use equations \eqref{eq18}
\begin{align*}
 c(x)&= c_{m,n}\int_0^x H(t)(1+t)t^{n/2}(1-t)^{-3} u_4(t)dt\\[2mm]
 d(x)&= -c_{m,n}\int_{\frac 12}^x H(t)(1+t)t^{n/2}(1-t)^{-3}u_3(t)dt +d_0\,.
\end{align*}
If $b\neq 0$, then $H(t) =a_0u_1(t)+bu_2(t)$ behaves like
$(1-t)^{2m+1-n}$ if $2m<n-1$ and like $\log(1-t)$ if $2m=n-1$,
resulting in $c(x)=O(1-x)^{2m-n-1}$, $d(x)=O(\log(1-x))$ if
$2m<n-1$ and $c(x)=O((1-x)^{-2}\log (1-x))$,
$d(x)=O((1-x)^{-1}\log(1-x))$ if $2m=n-1$. So if $b\neq 0$ one has
$G(x)=O(\log(1-x))$ if $2m<n-1$ and $G(x)=O((1-x)^{-1}\log (1-x))$
if $2m=n-1$. If $b=0$, then $H$ is bounded, giving
$c(x)=O((1-x)^{-2})$ and $d(x)=O(1)$ for $2m<n-1$,
$d(x)=O(\log(1-x))$ for $2m=n-1$. In case $2m<n-1$, however, we
can choose the constant $d_0$ so that $d(1)=0$, and then
$d(x)=O(1-x)^{n-2m-1}$. This choice gives $G(x)=O(1-x)^{n-2m-1}$
for $2m<n-1$. For $2m=n-1$, no choice of $d_0$ can improve the
bound $G(x)=O(\log(1-x))$.

It remains to estimate the growth of $A_2(r)$ near $r=1$. Recall
that the definition \eqref{eq14} of $A_3$ translates to
$A_2(r)=2G'(x)-(1-x)^{n-2m-2}H(x)$. Both terms growth like
$(1-x)^{n-2m-2}$, but a cancellation occurs. The functions $u_1,
u_3$ are $C^\infty$ at 1 and have developments
\begin{align*}
    u_3(x)&= A(1-x)^{n+1-2m} + O(1-x)^{n+2-2m}\\[2mm]
    u'_3(x)&= -A(n+1-2m)(1-x)^{n-2m} + O(1-x)^{n+1-2m}\\[2mm]
    H(x)&= a_0u_1(x)= B + O(1-x)\,.
\end{align*}
In $u_4(x)=F(m,m-\frac n2, 1+\frac n2, x)$, $1+\frac n2-m-(m-\frac
n2)=n+1-2m\geq 2$, whence $u_4$ has a finite derivative at 1 and a
development
\[
        u_4(x)= C+D(1-x)+O(1-x)^{1+\varepsilon}\qquad\forall
        \varepsilon<1\,, \quad u'_4(x)=O(1)\,.
\]
Then $W(x)= u'_3u_4-u_3u'_4=CA(2m-n-1)(1-x)^{n-2m}+\dotsc$, and so
the constant $c_{mn}$ in \eqref{eq18} is $CA(2m-n-1)$. Then from
\eqref{eq18}
\begin{align*}
c'(x)&= \frac{B(1-x)^{-3}}{A(2m-n-1)}+ O(1-x)^{-2}\\[2mm]
d'(x)&= -\frac{B(1-x)^{n-2m-2}}{C(2m-n-1)}+ O(1-x)^{n-2m-1}
\end{align*}
which gives
\begin{align*}
        c(x)&= \frac 12\frac{B}{2(2m-n-1)}(1-x)^{-2}+
                O(1-x)^{-1}\\[2mm]
        d(x)&= O(1-x)^{n-2m-1}\,, \,2m<n-1\,, \quad
        O(\log(1-x)), 2m=n-1\,.
\end{align*}
But $G'=c(x)u'_3(x)+d(x)u'_4(x)$; the second term $d(x)u'_4(x)$
satisfies the required bound while the first $c(x)u'_3(x)$ has a
development
\begin{align*}
   c(x)u'_3(x)&= -\frac 12 \frac{B}{A(2m-n-1)}
               A(n+1-2m)(1-x)^{n-2m-2}+ O(1-x)^{n-2m-1}\\[2mm]
   &=\frac{B}{2}(1-x)^{n-2m-2} + O(1-x)^{n-2m-1}\,.
\end{align*}
As $(1-x)^{n-2m-2}H(x)= B(1-x)^{n-2m-2}+ O(1-x)^{n-2m-1}$ the
bound for $A_2$ follows for $2m\leq n-1$.

However, for $m=\frac n2$, this no longer holds. Indeed, from
\eqref{eq17}, where $c=0$, $a=a_0$, $b-a$)
\[
        2G'(x)= x^{-\frac n2-1}\int_0^x s^{n/2}(1+s)(1-s)^{-3} a(s^{-\frac
        n2}-1)\,ds
\]
has development
\[
        2G'(x)= n a(1-x)^{-1}+ O(\log(1-x))
\]
while
\[
        (1-x)^{-2}a(x^{-\frac n2}-1)=\frac n2 a(1-x)^{-1}+\dotsc
\]

We point out that all this can be obtained, in loose terms,
working directly with the hypergeometric equations relating $G, H$
\[
        x(1-x)G''(x)+\left[\frac n2+1-(2m+1-\frac n2)x\right] G'(x)
        -m(m-\frac n2)G=\frac 12 (1+x)(1-x)^{n-2m-2}H(x)
\]
and using asymptotic developments. If $H(x)= h_0+h_1(1-x)+\dotsc$
and $G(x)=g_j(1-x)^j+\dotsc$, identifying the lower order terms in
both sides gives,
\[
        g_j j(j-1-n+2m)(1-x)^{j-1}=h_0(1-x)^{n-2m-2}.
\]
When $H\equiv 0$, one must have either $j=0$ (corresponding to
$u_4$) or $j=n-2m+1$ (corresponding to $u_3$). For the
inhomogeneous equation, if $j\neq 0$, $j\neq n+1+2m$ (that is $G$
contains no contribution from $u_3$, $u_4$) one finds $j=n-2m-1$
if $2m<n-1$ and $g_jj=-\frac {h_0}{2}$. Then $2G'(x)=
h_0(1-x)^{n-2m-2}+\dotsc$, $(1-x)^{n-2m-2}H(x)=h_0(1-x)^{n-2m-2}$,
showing cancellation. An analogous argument works if $2m=n-1$, but
not for $2m=n$.

We summarize the results in this and the previous subsections:

\begin{theorem}\label{th}
For $|n-2m|>1$, there is a unique doubly invariant kernel
$$k_m(x,y)= A_1(|\varphi_x y)\gamma_m + A_2(|\varphi_x y|)
\tau\wedge\gamma_{m-1}, m\neq 0; k_m(x,y)=
A(|\varphi_xy|)\gamma_m, m=0,n$$ for which \eqref{eq12} holds, and
satisfying moreover
\[
        |A_i(r)| =O(1-r^2)^{|n-2m|-1}\quad\text{as}\quad r\to 1\,.
\]
For $m=\frac{n\pm 1}{2}$, there is a one-parameter family of such
kernels satisfying
\[
        |A_i(r)|= O(\log (1-r^2))\,.
\]
For $m=\frac n2$, there is a one-parameter family of such kernels
satisfying
\[
        |A_i(r)|=O(1-r^2)^{-1}\,.
\]
In all cases $A_1(r)\sim r^{2-n}$, $A_2(r)\sim r^{-n}$ as $r\to
0$.
\end{theorem}

For $|n-2m|>1$, we call $k_m(x,y)$ the \emph{Riesz kernel for
$m$\forms} in $\Hn$, and
\[
        L\eta(x)=\int_{\Hn}\eta(y)\wedge *_yk_m(x,y)
\]
\emph{the Riesz potential of $\eta$}, whenever this is defined.
From \eqref{eq11} we see that
\begin{equation}\label{eq21}
 |k_m(x,y)| =O(1-r^2)^{n-m-1}.
\end{equation}
With the notations used before, the function $A_3(r)= H(r^2)$ is
bounded with bounded derivatives near $r=1$. Then \eqref{eq14} and
symmetry imply
\begin{equation}\label{eq22}
        |d_xk_m(x,y)|,\,|d_yk_m(x,y)|=O(1-r^2)^{n-m-1}
\end{equation}
too. The growth of $A_3$ also implies $A_4=O(1-r^2)^{n-2m}$
because $B_1\equiv 0$, and then \eqref{eq15} gives as well
\begin{equation}\label{eq23}
    |\delta_xk_m(x,y)|, \,|\delta_yk_m(x,y)| =O(1-r^2)^{n-m-1}\,.
\end{equation}

By construction, one has $L\Delta\eta=\eta$ for $\eta\in
\dcal_m(\Hn)$. We will need the following generalization of this
fact.

\begin{proposition}\label{pro}
If $\eta$ is a smooth form in $\Hn$ such that
\[
        |\eta(y)|,\,|\nabla\eta(y)| = o(1-|y|^2)^m\,,
        \quad y\in\Bn
\]
then $L\Delta\eta=\eta$.
\end{proposition}

\begin{proof}
In \eqref{eq19} we would get an extra term
\[
    \int_{|y|=R}(k_m\wedge \aster d\eta +\delta k_m\wedge\aster \eta -\delta\eta\wedge \aster k_m-\eta\wedge
    \aster dk_m)\,.
\]
Estimates \eqref{eq21}, \eqref{eq22} and \eqref{eq23} imply that
with $x$ fixed and $|y|=R\nearrow 1$
\[
        |k_m|,\,|\delta k_m|,\,|dk_m|=O(1-R^2)^{n-m-1}\,.
\]
Inserting $|\eta(y)|, |\nabla\eta(y)|=o(1-|y|^2)^m$ we see that
this extra term vanishes as $R\nearrow 1$.
\end{proof}

\section{Proof of the main theorem}\label{sec4}

\noindent {\bf 4.1.} Once the Riesz form $k_m(x,y)$ has been
found, our aim is now to prove that the corresponding convolution
\[
        L_m\eta(x)=\int_{\Hn}\eta(y)\wedge*_yk_m(x,y)
\]
satisfies
\begin{equation}\label{eq24}
        ||L_m\eta||_{p,s+2}\leq c||\eta||_{p,s}
\end{equation}
for $m\neq \frac{n\pm 1}{2}$, $\frac n2$, and $p$ in the range
$p_1(m)=\frac{n-1}{n-1-m}<p<\frac{n-1}{m}=p_2(m)$, and for a
compactly supported $m$\form\ $\eta$ (recall that we are assuming
without loss of generality that $m\leq\frac n2$). Since these are
dense in the Sobolev spaces and we already know that $\Delta
L_m\eta=L_m\Delta\eta=\eta$, this will prove the theorem for
$m\neq\frac{n\pm 1}{2}$, $\frac n2$. The case $m=\frac{n\pm 1}{2}$
will be commented later.

We work in the translation invariant basis $w^I$. Taking into
account formulas \eqref{eq9} and (10) for $\gamma$, $\tau$, the
Riesz- form is written in the $\Rbb_+^n$ model
\[
        k_m(x,y)=\sum_{|I|=|J|=m} a_{I,J}(S_y x) w^I(x)\otimes w^J(y)
\]
where each coefficient $a_{I,J}$ has an expression, with $z=S_yx$
\[
        a_{I,J}(z)=\Psi_{I,J}(r)\frac{p_{I,J}(z)}{(|z|^2+2z_n+1)^{2m}}\,,
        \quad
        r^2=\frac{1+|z|^2-2z_n}{1+|z|^2+2z_n} =
        \frac{|x-y|^2}{|x-y|^2+4x_ny_n}\,.
\]
Here $p_{I,J}(z)$ is a certain polynomial in $z_1,\dots,z_n$,
$\Psi_{I,J}$ is $C^\infty$ in $(0,1)$ with $\Psi_{I,J}(r)\sim c_0
r^{2-n}$ as $r\searrow 0$, $\Psi_{I,J}(r)=O(1-r^2)^{n-m-1}$ as
$r\nearrow 1$. The term $q_{I,J}(z)=p_{I,J}(z)/(|z|^2 +2 z_n
+1)^{2m}$ is bounded.

If $\eta=\sum_I \eta_I(y)w^I(y)$, the coefficient $(L\eta)_I(x)$
of $L\eta$ in the basis $w^I$ is a finite linear combination of
hyperbolic convolutions
\[
        (L\eta)_I(x)=\sum_J\int_{\Hn} \Psi_{I,J}(r) q_{I,J}(z)
        \eta_J(y)\,d\mu(y)\,.
\]
By ellipticity of $\Delta$, $L\eta$ is a smooth form. Moreover,
since $\eta$ has compact support, we see from \eqref{eq2} and
\eqref{eq21}, \eqref{eq22}, \eqref{eq23} that, in the ball model,
\[
        |L\eta(x)|,\,|d(L\eta)(x)|,\,|\delta(L\eta)(x)|=O(1-|x|^2)^{n-m-1}
\]
which amounts to
\begin{equation}\label{eq25}
        |(L\eta)_I(x)|,\,|X_i(L\eta)_I(x)|=O(1-|x|^2)^{n-m-1}\,.
\end{equation}
We claim that for second-order derivatives we have too
\begin{equation}\label{eq26}
    |X_jX_i(L\eta)_I(x)|=O(1-|x|^2)^{n-m-1}\,,\quad i.e.
    \quad|\nabla^{(2)}(L\eta)(x)|=O(1-|x|^2)^{n-m-1}\,.
\end{equation}
Notice that since we already know that $\Delta L\eta=\eta$, from
the expression of $\Delta$ in the basis $w^I$ given in (3) it
follows that it is enough to show that for $j<n$. We will see
below (equation (30) and invariance of the $X_i$) that each of the
functions $a(z)=\Psi_{I,J}(r)q_{I,J}(z)$ satisfies
\[
 |X_jX_i a(z)| =O(1-r^2)^{n-m-1}
\]
from which \eqref{eq26} follows as before. In fact, the discussion
that follows will show that
$|\nabla^{(k)}L\eta(x)|=O(1-|x|^2)^{n-m-1}$ $\forall k$.

 We continue the proof of \eqref{eq24}. We claim first that it is enough to prove \eqref{eq24} for $s=0$.
For a smooth form $\eta=\sum\eta_I w^I$ let $X_i\eta$ denote here
the $m$\form\ $X_i\eta=\sum X_i\eta_Iw^I$. It is clear from
formulas (3) and the commutation properties
\[
        [X_i,X_j]=0\,,\qquad i,j<n\,,\qquad [X_1,X_i]=X_i\,,\qquad
        i<n
\]
that for each $i$ there is an operator $P_i$ of order two in the
$X_1,\dots, X_n$ such that
\[
        X_i\Delta\eta-\Delta(X_i\eta)=P_i(X)\eta\,.
\]
Applying this to $L\eta$, which is smooth by the ellipticity of
$\Delta$, we get $(X_i-\Delta X_i L)\eta=P_i(X)L\eta$. But
$X_iL\eta$ satisfies, by \eqref{eq25} and \eqref{eq26}
\[
        |X_iL\eta(x)|,\,|d(X_iL\eta)(x)|,\,|\delta(X_iL\eta)(x)|=O(1-|x|^2)^{n-m-1}
\]
and hence by Proposition 3.2, $L\Delta =\Id$ on it. We conclude
that for all $\eta\in\dcal_m(\Hn)$
\[
        (LX_i-X_iL)\eta=LP_i(X)L\eta\,.
\]
Assume that \eqref{eq24} has been proved up to $s$, so that by
density it holds for $\alpha\in H^s_{m,p}(\Hn)$ too, and let
$\gamma$ be a multiindex of length $|\gamma|\leq s$. For
$i=1,\dots,n$ and $\eta\in\dcal_m(\Hn)$,
\[
        X^\gamma X_iL\eta=X^\gamma LX_i\eta-X^\gamma L P_i(X)L\eta
\]
so using twice the induction hypothesis
\begin{align*}
        ||X^\gamma X_iL\eta||_p&\leq\const\left(||X_i\eta||_{p,s} +
        ||P_i(X)L\eta||_{p,s}\right)\leq\\[2mm]
        &\leq \const\left(||\eta||_{p,s+1} + ||\eta||_{p,s}\right)
\end{align*}
proving \eqref{eq24} for $s+1$. Proving \eqref{eq24} for $s>0$
means proving
\[
        ||(L\eta)_I||_p\,,||X_i(L\eta)_I||_p\,,||X_jX_i(L\eta)_I||_p\leq
        \const ||\eta||_p\,.
\]
As before, using that we already know that $\Delta L\eta=\eta$ we
see that for the second-order derivatives we may assume $j<n$. In
the following we delete the indexes $I,J$ and denote by
$a(z)=\psi(r)Q(z)$ a convolution kernel with $\psi, Q$ as above,
and proceed to prove that the convolution
\[
        (C_a\alpha)(z)=\int_{\Hn} a(S_yx)\alpha(y)\,d\mu(y)
\]
satisfies
\begin{equation}\label{eq27}
        ||C_a\alpha||_p, ||X_i(C_a\alpha)||_p,
        ||X_jX_iC_a(\alpha)||_p\leq\const ||\alpha||_p\,,\qquad
        p_1\leq p\leq p_2
\end{equation}
where in the last case we may assume that $j<n$. The fields $X_i$
are invariant, and therefore $X_i C_a\alpha$, $X_jX_iC_a\alpha$
are obtained, respectively, by convolution with $Z_ia$, $Z_jZ_i a$
(by \eqref{eq4}). Recall that
$\psi(r)=O(1-r^2)^{n-m-1}=O\left(\frac{4z_n}{1+|z|^2+2z_n}\right)^{n-m-1}$
as $r\nearrow 1$ and $\psi(r)\sim r^{2-n}$ as $r\searrow 0$.

In order to estimate $Z_ia, Z_iZ_ja$, we collect first some
auxiliary estimates. We claim that
\begin{align}\label{eq28}
   & |Z_i Q|\leq \const\,,\qquad |Z_iZ_j Q|\leq \const\nonumber\\
   &\ \\[-2mm]
   & |Z_ir|\leq\const (1-r^2)\,,\qquad |Z_iZ_j r|\leq \const
   r^{-1}(1-r^2)\,.\nonumber
\end{align}
The first two are routinely checked, for instance, when
differentiating the denominator in $Q$,
\[
        \left|Z_i\frac{1}{(1+|z|^2+2z_n)^{2m}}\right|
        =\left|\frac{4mz_nz_1}{(1+|z|^2+2z_n)^{2m+1}}\right|
        \leq \frac{\const}{(1+|z|^2+2z_n)^{2m}} \qquad (i<n)
\]
so that the term $p_{I,J}(z)Z_i\left[(1+|z|^2+2z_n)^{-2m}\right]$
will still be bounded. All other terms can be treated similarly.
Differentiating $1-r^2=\frac{4z_n}{1+|z|^2+2z_n}$ we get
\begin{align*}
 Z_ir&= \frac{1-r^2}{2}\,\frac{z_iz_n}{r(1+|z|^2+2z_n)}\,,\qquad
 Z_nr = -\frac{1-r^2}2 \frac 1r
\frac{1+|z|^2-2z_n^2}{1+|z|^2+2z_n}\\[4mm]
Z_jZ_ir &=\frac{1-r^2}{2r}\left\{
        \frac{\delta_{ij}z_n^2}{1+|z|^2+2z_n}- \frac{1+5r^2}{2r^2}\,
        \frac{z_iz_jz_n^2}{(1+|z|^2+2z_n)^2}\right\}\,,
        \qquad i,j<n\\[4mm]
Z_jZ_nr&=\frac{1-r^2}{r}\left\{-\frac{2z_n^2z_j(1+z_n)}{(1+|z|^2+2z_n)^2)}
    +\frac{(1+r^2)}{4r^2}\,\frac{z_jz_n(1+|z|^2-2z_n)}{(1+|z|^2+2z_n)^2}\right\}\,,\qquad
    j<n\,.
\end{align*}
These imply \eqref{eq28} because
\[
        |z_iz_n|, 1+|z|^2-2z_n^2\leq
(1+|z|^2-2z_n)^{1/2}(1+|z|+2z_n)^{1/2}=r(1+|z|^2+2z_n)\,.
\]
Now
\begin{align}\label{eq29}
    Z_i a(z)&= \psi'(r) Z_ir \,Q(z)+ \psi(r)(Z_i Q)(z)\nonumber\\
\ \\[-3mm]
    Z_jZ_ia(z)&=\psi''(r)(Z_ir)(Z_jr)Q(z) + \psi'(r)(Z_jZ_ir)Q+
    \psi'(r)Z_irZ_j Q\nonumber \\[2mm]
    &\qquad + \psi'(r)Z_jr Z_iQ + \psi(r)Z_jZ_i Q\,.\nonumber
\end{align}
The estimates \eqref{eq28} imply
\begin{align}\label{eq30}
   & |a(z)|\,,|Z_ia(z)|\,,
  |Z_jZ_ia(z)|=O(1-r^2)^{n-m-1}\qquad\text{as}\quad r\nearrow 1\nonumber\\
   \ \\[-3mm]
  &|a(z)| = O(r^{2-n})\,, |Z_ia(z)|=O(r^{1-n})\,, |Z_jZ_i
  a(z)|=O(r^{-n})\nonumber
\end{align}

We will call a convolution kernel $b(z)$ \emph{$m$\admissible\ }
if $|b(z)|=O(r^{1-n})$ as $r\searrow 0$ and, moreover,
$|b(z)|=O(1-r^2)^{n-m-1}$ as $r\nearrow 1$. We will prove later
(Theorem 4.2.) that hyperbolic convolution with $m$-admissible
kernels defines a bounded operator in $L^p(\Hn)$ for the range
$(p_1(m), p_2(m))$ specified in the statement of the main result.
From the estimates \eqref{eq29} we see that $a$ and $Z_ia$ are
$m$\admissible\ kernels, and so \eqref{eq27} will be proved for
them. As $|Z_jZ_ia(z)| =O(r^{-n})$ has the critical non-integrable
singularity at $r=0$, $Z_j Z_ia(z)$ is not an $m$\admissible
kernel. Notice however from \eqref{eq29}, \eqref{eq30}
 that the last three terms $\psi'(r)Z_ir Z_jQ$, $\psi'(r)Z_jrZ_i
 Q$, $\psi(r)Z_jZ_iQ$ are indeed $m$\admissible.
Moreover, the estimate $|Z_i Q|\leq \const$ implies that $Q$ is
Lipschitz with respect the hyperbolic metric, in particular
\[
        Q(z)=Q(e)+ O(\log\frac{1+r}{1-r}) =Q(e)+O(r)
\]
for small $r$. This means that replacing $Q$ by $Q-Q(e)$ in the
first two terms leads to an $m$\admissible\ kernel again. All this
leaves us with the kernel
\[
        \psi''(r)Z_irZ_jr + \psi'(r)Z_jZ_ir \,,\qquad j<n\,.
\]
If $\psi(r)=c_0r^{2-n}+\dotsc$, write
$\phi(r)=c_0r^{2-n}(1-r^2)^{n-m-1}$; then the above differs from
\[
        \phi''(r)Z_irZ_jr + \phi'(r)Z_jZ_ir
\]
in an $m$\admissible\ kernel. By the same reason, we may replace
$\phi''(r)$, $\phi'(r)$ respectively by
$(r^{2-n})''(1-r^2)^{n-m-1}$, $(r^{2-n})'(1-r^2)^{n-m-1}$, that is
to say we must deal with the convolution kernel
\begin{equation}\label{eq31}
    (1-r^2)^{n-m-1} Z_jZ_i(r^{2-n})\,.
\end{equation}

We introduce a class of singular hyperbolic convolution kernels to
deal with the later. For this purpose it is more convenient to
work in the ball model, so now $b$ is defined in $\Bn$ and
$r=|z|$. We replace the integrable singularity $r^{1-n}$ by a
typical Calder\'{o}n-Zygmund singularity (see e.g. \cite{S}). Thus, we
will call $b$ a \emph{$m$-Calder\'{o}n-Zygmund singular kernel} if it
has the form
\[
    b(z)=\Omega(w)r^{-n}(1-r^2)^{n-m-1}\,,\quad z=rw\,, \quad w\in S^{n-1}
\]
where $\Omega$ is say a Lipschitz function on $S^{n-1}$ satisfying
the cancellation condition
\begin{equation}\label{eq32}
    \int_{S^{n-1}}\Omega(w)\,d\sigma(w)=O\,.
\end{equation}
In Theorem 4.2. below we prove that $m$-Calder\'{o}n-Zygmund singular
kernels define bounded operators in the same range of $p$. With
the following proposition, applied to $\phi_2(z)=|z|^{2-n}$ this
will end the proof of the main result. The proposition is the
analogue of the well-known statement that for $\phi$ smooth and
homogeneous of degree $1-n$ in $\Rbb^n$,
$\frac{\partial\phi}{\partial x_i}$ defines a Calder\'{o}n-Zygmund
kernel (it is homogeneous of degree $-n$ and the cancellation
condition (32) is automatically satisfied, because
$\int_{r_1<|x|<r_2}\frac{\partial\phi}{\partial
x_i}\,dV(x)=\left(\int_{|x|=r_2}-\int_{|x|=r_1}\right)\phi(x)\,dx^1\wedge\dots\wedge
dx^i\wedge\dots\wedge dx^n=0$).

\begin{proposition}\label{pro4.1.}
If $\phi_1,\phi_2$ are homogeneous functions of degree $1-n$,
$2-n$ respectively, the kernels $(1-r^2)^{n-m-1}Z_i\phi_1$,
$(1-r^2)^{n-m-1}Z_jZ_i\phi_2$ are sum of $(m-1)$\admissible\ and
$(m-1)$-Calder\'{o}n-Zygmund singular kernels.
\end{proposition}

\begin{proof}
We replace the $Z_j$ by $Y_j=(1-r^2)\frac{\partial}{\partial
z_j}$; we have
\begin{align*}
  Y_i(\phi_1)&= (1-r^2)\frac{\partial\phi_1}{\partial
  z_i}\qquad\quad Y_i\phi_2=(1-r^2)\frac{\partial \phi_2}{\partial
  z_i}=(1-r^2)O(r^{1-n})\\[2mm]
    Y_jY_i\phi_2 &= (1-r^2)\frac{\partial^2\phi_2}{\partial z_i\partial z_j}-2(1-r^2) z_j
    \frac{\partial\phi_2}{\partial z_i}=(1-r^2)\frac{\partial^2\phi_2}{\partial z_i\partial z_j}
    +(1-r^2)O(r^{2-n})
\end{align*}
so in all cases we get an extra factor $(1-r^2)$. Besides
$\frac{\partial\phi_1}{\partial z_i}$,
$\frac{\partial^2\phi_2}{\partial z_i\partial z_j}$ are, as noted
before, homogeneous of degree $-n$, and satisfy the cancellation
condition \eqref{eq32}.
\end{proof}

\vspace{1cm}

\noindent {\bf 4.2.} It remains to prove
\begin{theorem}\label{theo4.2}
Both $m$-admissible  and $m$-Calder\'{o}n-Zygmund kernels define by
hyperbolic convolution bounded operators in $L^p(\Hn)$ for
$\frac{n-1}{n-1-m}<p<\frac{n-1}{m}$, $0\leq m<\frac{n-1}{2}$.
\end{theorem}

We will make use of the following well-known Schur's lemma for
boundedness in $L^p$ of an integral operator with positive kernel.

\begin{lemma}\label{le4.3}
If $K(x,y)$ is a positive kernel in a measure space $X$ and
$1<p<\infty$, the operator $Kf(x)=\int_XK(x,y)f(y)\,d\mu(y)$ is
bounded in $L^p(\mu)$ if and only if there exists $h\geq 0$ such
that
\begin{align}
    \int_XK(x,y)h(y)^{q}\,d\mu(y)&=O(h(x)^q)\,,\qquad x\in X \label{eq33}\\[3mm]
    \int_XK(x,y)h(x)^p \,d\mu(x)&= O(h(y)^p)\,,\qquad y\in
    Y\label{eq34}\,.
\end{align}
Here $q$ is the conjugate exponent of $p$, $\frac 1p+\frac 1q=1$.
If $h$ can be taken $\equiv 1$ that is
\[
        \sup_x\int_XK(x,y)\,d\mu(y)\,,\sup_y\int_XK(x,y)\,d\mu(x)<+\infty
\]
then $K$ is bounded in $L^p(\mu)$ for all $p$, $1\leq
p\leq\infty$.
\end{lemma}

\begin{proof} Let us prove Theorem 4.2. If $b$ is $m$\admissible, $b=b_1+b_2$ with $b_1(z)=O(r^{1-n})$ for
$r\leq 1/2$, $b_1(z)=0$ for $r>1/2$, and $b_2(z)=O(1-r^2)^{n-m-1}$
for all $r$. We apply to $b_1$ the second criteria in
Lemma~\ref{le4.3} working in the ball model (recall that $|S_y
X|=|\varphi_yx|$ is symmetric in $x,y$)
\begin{align*}
 &\int_Xb_1(S_yx) \,d\mu(x)\,,\int_X b_1(S_y x)\,d\mu(y)\leq
 c\int_{|S_yx|\leq 1/2}|S_yx|^{1-n}d\mu(x) \\[2mm]
 &=c\int_{|z|\leq1/2}|z|^{1-n}d\mu(z)=\const\int_0^{1/2}\frac{dr}{(1-r^2)^n}<+\infty.
\end{align*}
We apply to $(1-r^2)^{n-m-1}$ the criteria of the first part on
Lemma~\ref{le4.3}, working this time for convenience in the
half-space model, where the kernel is written
\[
        K(x,y)=(1-r^2)^{n-m-1}=\left(\frac{4z_n}{1+|z|^2+2z_n}\right)^{n-m-1}
        =\left(\frac{4x_n y_n}{|x-y|^2 + 4x_n
        y_n}\right)^{n-m_1}\,.
\]
We test $h(y)=y_n^\alpha$ in \eqref{eq33} for an exponent $\alpha$
to be chosen, so we need
\[
        \int_{y_n>0}\frac{y_n^{-k-1+\alpha q}\,dy}{(|x-y|^2 +
        4x_ny_n)^{n-k-1}}= O\left(x_n^{\alpha q+k+1-n}\right)\,.
\]
We write $|x-y|^2+4x_ny_n=|x'-y'|^2 + (x_n+y_n)^2$, where
$x'=(x_1,\dots,x_{n-1})\in\Rbb^{n-1}$ and analogously for $y'$,
and integrate first in $y'$. One has for $2m<n-1$
\begin{align*}
        \int_{\Rbb^{n-1}}\frac{dy'}{(|x'-y'|^2 + (x_n+y_n)^2)^{n-m-1}}
    &=c\int_0^\infty\frac{s^{n-2}}{(s^2+(x_n+y_n)^2)^{n-m-1}}\\[3mm]
    &= O((x_n+y_n)^{2m+1-n})
\end{align*}
and so the above becomes
\[
        \int_0^\infty\frac{y_n^{\alpha
        q-m-1}dy_n}{(x_n+y_n)^{n-1-2m}}= O(x_n^{\alpha q+m+1-n})\,.
\]
By homogeneity ($y_n=x_nt$) this reduces to
\[
        \int_0^\infty\frac{t^{\alpha q-m-1}}{(1+t)^{n-1-2m}}=O(1)
\]
which holds whenever $m<\alpha q<n-1-m$. By symmetry, for
\eqref{eq34} we need as well $m<\alpha p<n-1-m$. Therefore, a
choice of $\alpha$ is possible whenever $m\max\left(\frac 1p,\frac
1q \right)< (n-1-n)\min\left(\frac 1p,\frac 1q\right)$ and this
gives the range $\frac{n-1}{n-1-m}<p<\frac{n-1}{m}$.

Consider now a $m$-Calder\'{o}n-Zygmund kernel
$b(z)=\Omega(w)r^{-n}(1-r^2)^{n-m-1}$. Since
$|S_yx|=|\varphi_yx|$, we may replace $z=S_yx$ by $z=\varphi_x y$.
Using \eqref{eq1} this is given by
\[
    z=\frac{(x-y)(1-|x|^2)+ x|x-y|^2}{A}
\]
where we use the notation $A=(1-|x|^2)(1-|y|^2)+|x-y|^2$; note
that
$$(1-|x|^2),
(1-|y|^2)\lesssim A^{1/2}$$
Also recall that $r=|z|$ and
$Ar^2=|x-y|^2$. Hence we can write
\[
         \frac z r-\frac{x-y}{|x-y|}=\frac{x-y}{|x-y|}\left(\frac{1-|x|^2}{\sqrt{A}}-1\right)
         +x\cdot r\,.
\]
But
\[
    \frac{1-|x|^2}{\sqrt{A}}-1 =\frac{(1-|x|^2)^2-A}{\sqrt{A}((1-|x|^2)+\sqrt{A})}
    =\frac{(1-|x|^2)O(|x-y|)+O(|x-y|^2)}{A}
\]
is $O(r)$. Therefore, modulo an $m$\admissible\ kernel, we may
replace $\Omega(w)$ by $\Omega\left(\frac{x-y}{|x-y|}\right)$.
This leaves us with the kernel
\[
    K=(1-r^2)^{n-m-1}\Omega\left(\frac{x-y}{|x-y|}\right)r^{-n}
    =(1-r^2)^{n-m-1}|x-y|^{-n}\, \Omega\left(\frac{x-y}{|x-y|}\right)
    A^{\frac n2}(x,y)\,.
\]
Fix $p$, $1<p<\infty$. Write
\[
    A^{\frac n2}(x,y)= (1-|x|^2)^{\frac np}(1-|y|^2)^{\frac nq} +
    O\left(|x-y|\, A^{\frac{n-1}2}\right)
\]
Since $|x-y|^{1-n} A^{\frac{n-1}2}=r^{1-n}$, the kernel $K$
differs from
\[
        (1-r^2)^{n-m-1}|x-y|^{-n}\Omega\left(\frac{x-y}{|x-y|}\right)
        (1-|x|^2)^{\frac np}(1-|y|^2)^\frac nq
\]
in a $m$\admissible\ kernel, so we keep this one. We write it as
the sum of
\[
    |x-y|^{-n}\Omega\left(\frac{x-y}{|x-y|}\right)
    (1-|x|^2)^{\frac np} (1-|y|^2)^{\frac nq}=K_1 (x,y)
\]
and another $K_2(x,y)$ which we estimate by
\[
        |K_2(x,y)| = O\left( r^2|x-y|^{-n}(1-|x|^2)^{\frac np}(1-|y|^2)^{\frac nq}\right)
        = O\left(r^{2-n}(1-|x|^2)^{\frac np}(1-|y|^2)^{\frac nq}A^{-\frac
        n2}\right)\,.
\]
Write $K_\Omega$ for the (euclidean) Calder\'{o}n-Zygmund
convolution operator with kernel \linebreak
$|x-y|^{-n}\Omega\left(\frac{x-y}{|x-y|}\right)$, which as it is
well-known, satisfies an $L^p(dV)$-estimate. Notice that
\[
        K_1f(x)= (1-|x|^2)^{\frac np} K_\Omega\left(f(1-|y|^2)^{-\frac np}\right)
\]
and therefore, using the $L^p$-boundedness of $K_\Omega$
\[
        \int_{\Bn}|K_1f(x)|^p \,d\mu(x)=\int_{\Bn}
        \left|K_\Omega(f(1-|y|^2)^{-\frac np})\right|^p
        \,dV(x)\leq \int_{\Bn}|f(x)|^p\,d\mu(y)\,.
\]
For $K_2$, we can ignore the integrable simgularity $r^{2-n}$ and
arguing as we just did with $K_1$, we need to show that the
integral operator
\[
        K_3 f(x)=\int_{|y|\leq 1}\frac{1}{(1-|x|+|x-y|)^n}
        f(y)\,dV(y)
\]
satisfies $L^p(dV)$-estimates for all $p$, $1<p<\infty$. To see
this, just check that the criteria in Lemma~\ref{le4.3} holds with
$h(x)=(1-|x|^2)^{-\frac 1{pq}}$

\end{proof}

Notice that in case $m=0$ a $m$-Calder\'{o}n-Zygmund kernel defines a
bounded operator in all $L^p(\Hn)$, $1<p<\infty$: this is the
right analogue of the euclidian  kernels, because $(1-r^2)^{n-1}$
is the typical growth at infinity of a weak $L^1(d\mu)$ function
in $\Hn$.

\vspace{1cm}

\noindent {\bf 4.3.} Finally we make some comments, with no
proofs, on the critical case $m=\frac{n-1}{2}$ in the main
Theorem. In this case the $m$\admissible\ and $m$-Calder\'{o}n-Zygmund
operators appearing in $X_jX_i C_au$, etc. have
$(1-r^2)^{\frac{n-1}2}\log\frac{1}{1-r^2}$ instead of
$(1-r^2)^{n-m-1}=(1-r^2)^{\frac{n-1}2}$ as a factor. One can then
prove that for $\beta>0$ and $2\leq p <2+\frac{2\beta}{(n-1)}$,
\[
        ||L_p\eta||_{p,2}\leq \const\int_{\Bn}|\eta|^p
    (1-|y|^2)^{-\beta}\,d\mu(y)\,.
\]
The $L^p$-estimates do not hold in this case for any $p$, because
they do not hold for $p=2$  and $\Delta$ is self-adjoint.

\vspace{1cm} {\small \noindent
Departament de Matem\`{a}tiques\\
Universitat Aut\`{o}noma de Barcelona\\
Campus de Bellaterra\\
08193 Cerdanyola del Vall\`{e}s (Barcelona), Spain\\
\textsf{bruna@mat.uab.es}}
\end{document}